\documentclass[11pt]{article}
\usepackage{amssymb,amsthm,amsmath,amscd}

\newtheorem{thm}{Theorem}[section]
\newtheorem{lem}[thm]{Lemma}

\renewcommand{\phi}{\varphi}

\newcommand{\Homeo}{\operatorname{Homeo}}
\newcommand{\Coinv}{\operatorname{Coinv}}
\newcommand{\id}{\operatorname{id}}

\newcommand{\N}{\mathbb{N}}
\newcommand{\Z}{\mathbb{Z}}
\newcommand{\Q}{\mathbb{Q}}
\newcommand{\R}{\mathbb{R}}
\newcommand{\T}{\mathbb{T}}

\newcommand{\xp}{(X,\phi)}
\newcommand{\yp}{(Y,\psi)}

\title{Torsion in coinvariants of certain Cantor minimal $\Z^2$-systems}
\author{Hiroki Matui 
\thanks{Supported in part by a grant 
from the Japan Society for the Promotion of Science} \\
Graduate School of Science \\
Chiba University \\
1-33 Yayoi-cho, Inage-ku, Chiba 263-8522, Japan}
\date{}

\begin{document}
\maketitle

\begin{abstract}
Let $G$ be a finite abelian group. 
We will consider a skew product extension of 
a product of two Cantor minimal $\Z$-systems 
associated with a $G$-valued cocycle. 
When $G$ is non-cyclic and the cocycle is non-degenerate, 
it will be shown that the skew product system has torsion 
in its coinvariants. 
\end{abstract}

\section{Introduction}

Let $X$ be a Cantor set and 
let $\Gamma$ be a countable discrete abelian group which acts on $X$. 
We write the action of $\Gamma$ by $\phi:\Gamma\to\Homeo(X)$. 
We define the group of coinvariants of $\xp$ by 
\[ \Coinv\xp=C(X,\Z)/
\langle\{f-f\circ\phi^{-\gamma}
:f\in C(X,\Z),\gamma\in\Gamma\}\rangle, \]
where $C(X,\Z)$ denotes 
the set of integer-valued continuous functions on $X$. 
When $G$ is $\Z^d$, 
the group of coinvariants is known to be isomorphic to 
the top-dimensional dynamical cohomology $H^d(\Z^d,C(X,\Z))$ of $\xp$. 
The aim of this paper is to compute the torsion part of 
$\Coinv\xp$ for a certain minimal action of $\Z^2$. 

G\"ahler \cite{G} discovered unexpected 5-torsion 
in the top-dimensional cohomology of a popular two-dimensional tiling 
which is called the T\"ubingen Triangle Tiling (TTT). 
The TTT is one of canonical projection tilings. 
For these tilings, Forrest, Hunton and Kellendonk \cite{FHK} related 
the cohomology of tilings to a certain group homology, 
and they expressed its rational rank in terms of the data defining the tiling 
in the cases of codimension at most three. 
G\"ahler, Hunton and Kellendonk \cite{GHK} extended this result and 
computed the torsion in the cohomology 
for such tilings of codimension at most three. 
In this paper, we would like to show that 
torsion in the top cohomology can arise in a much elementary setting. 
In addition, we will see that 
an example from \cite{GHK} with torsion can also be described 
in this context. 
By our method, however, 
we cannot `explain' why the 5-torsion arises in the cohomology of the TTT. 
It is worth clarifying its reason in terms of dynamical systems. 

Our setting is as follows. 
Let $X$ and $Y$ be Cantor sets and 
let $\phi\in\Homeo(X)$ and $\psi\in\Homeo(Y)$ be minimal homeomorphisms. 
Then $\phi\times\id$ and $\id\times\psi$ induce 
a minimal $\Z^2$-action on $X\times Y$. 
This is called the product of $\xp$ and $\yp$, 
and it is well-known that the group of coinvariants is isomorphic to 
the tensor product of $\Coinv\xp$ and $\Coinv\yp$. 
In particular, it is torsion-free. 
In this article, however, 
we will see that a skew product extension of the product system can 
have torsion in its coinvariants. 
More precisely, 
when $G$ is a non-cyclic finite abelian group, a skew product extension 
associated with a non-degenerate $G$-valued cocycle always has torsion 
in its coinvariants. 
Indeed, the torsion part of coinvariants depends only on 
the finite abelian group $G$.

\section{Torsion of a certain abelian group}

Let $G$ be a finite abelian group. 
We use the multiplicative notation for the group operation in $G$. 
Let $A$ and $B$ be finite sets and 
let $\mu:A\cup B\to G$ be a map. 
Suppose that both $\mu(A)$ and $\mu(B)$ generate $G$. 
Put 
\[ {\cal M}(A,B)=\Z[G]\otimes\Z^A\otimes\Z^B. \]
We denote each element of ${\cal M}(A,B)$ by 
\[ r=\sum_{a\in A}\sum_{b\in B}
r(a,b)\otimes a\otimes b\in{\cal M}(A,B), \]
where $r(a,b)$ is an element of the group ring $\Z[G]$. 

Let ${\cal A}\subset{\cal M}(A,B)$ be 
the subgroup consisting of those elements 
\[ \sum_{a\in A}\sum_{b\in B}
\alpha(a)(e-\mu(b))\otimes a\otimes b\in{\cal M}(A,B) \]
with $\alpha:A\to\Z[G]$. 
In a similar fashion, 
let ${\cal B}\subset{\cal M}(A,B)$ be 
the subgroup consisting of those elements 
\[ \sum_{a\in A}\sum_{b\in B}
\beta(b)(e-\mu(a))\otimes a\otimes b\in{\cal M}(A,B) \]
with $\beta:B\to\Z[G]$. 
Put 
\[ {\cal N}(A,B)={\cal M}(A,B)/({\cal A}+{\cal B}) \]
and let $\rho_{A,B}:{\cal M}(A,B)\to{\cal N}(A,B)$ be the quotient map. 
We denote by $T({\cal N}(A,B))$ the torsion part of ${\cal N}(A,B)$ 
as a $\Z$-module, that is, 
\[ T({\cal N}(A,B))=\{r\in{\cal N}(A,B):nr=0\text{ for some }n\in\N\}. \]
The aim of this section is to compute this group. 

Let ${\cal M}_\Q(A,B)=\Q[G]\otimes\Z^A\otimes\Z^B$. 
We regard ${\cal M}(A,B)$ as a subgroup of ${\cal M}_\Q(A,B)$. 
By replacing $\Z[G]$ with $\Q[G]$, 
we can similarly define ${\cal A}_\Q$ and ${\cal B}_\Q$ 
as subgroups of ${\cal M}_\Q(A,B)$. 
It is not hard to see 
\[ \rho_{A,B}^{-1}(T({\cal N}(A,B)))=\{s+t\in{\cal M}(A,B):
s\in{\cal A}_\Q\text{ and }t\in{\cal B}_\Q\}. \]

\begin{lem}\label{ABpure}
We have ${\cal A}_\Q\cap{\cal M}(A,B)={\cal A}$ and 
${\cal B}_\Q\cap{\cal M}(A,B)={\cal B}$. 
\end{lem}
\begin{proof}
It suffices to show only ${\cal A}_\Q\cap{\cal M}(A,B)={\cal A}$. 
It is obvious that 
${\cal A}$ is contained in ${\cal A}_\Q\cap{\cal M}(A,B)$. 
Let $\alpha:A\to\Q[G]$ be a map and 
suppose that 
\[ r=\sum_{a\in A}\sum_{b\in B}
\alpha(a)(e-\mu(b))\otimes a\otimes b \]
belongs to ${\cal M}(A,B)$. 
We observe $\alpha(a)(e-\mu(b))\in\Z[G]$ for all $a\in A$ and $b\in B$. 
Since $\mu(B)$ generates $G$, for each $a\in A$, 
there exists $s_a\in\Q$ such that 
\[ \alpha(a)-s_aN\in\Z[G], \]
where $N$ is $\sum_{g\in G}g\in\Z[G]$. 
By replacing $\alpha(a)$ with $\alpha(a)-s_aN$, 
we can see that $r$ belongs to ${\cal A}$. 
\end{proof}

At first, we would like to show that 
$T({\cal N}(A,B))$ does not depend on the choice of 
$A$, $B$ and $\mu:A\cup B\to G$. 
Let $C$ and $D$ be other finite sets and 
let $\nu:C\cup D\to G$ be a map. 
In exactly the same way, 
we can define ${\cal M}(C,D)$, ${\cal C}$, ${\cal D}$ and so on. 
Since $\mu(A)$ generates $G$, we can find 
$\{s(a,c)\}_{a\in A,c\in C}\subset\Z[G]$ such that 
\[ \sum_{a\in A}(e-\mu(a))s(a,c)=e-\nu(c) \]
for all $c\in C$. 
Similarly, since $\mu(B)$ generates $G$, we can find 
$\{t(b,d)\}_{b\in B,d\in D}\subset\Z[G]$ such that 
\[ \sum_{b\in B}(e-\mu(b))t(b,d)=e-\nu(d) \]
for all $d\in D$. 
We can define a homomorphism $\pi$ 
from ${\cal M}_\Q(A,B)$ to ${\cal M}_\Q(C,D)$ by 
\[ \pi(r\otimes a\otimes b)
=\sum_{c\in C}\sum_{d\in D}rs(a,c)t(b,d)\otimes c\otimes d \]
for $r\in\Q[G]$ and $a\in A$, $b\in B$. 
Clearly $\pi$ carries ${\cal M}(A,B)$ to ${\cal M}(C,D)$. 
It is also easily checked that 
$\pi({\cal A})\subset{\cal C}$, $\pi({\cal B})\subset{\cal D}$, 
$\pi({\cal A}_\Q)\subset{\cal C}_\Q$ and 
$\pi({\cal B}_\Q)\subset{\cal D}_\Q$. 
In particular, we obtain a homomorphism $\tilde\pi$ 
from ${\cal N}(A,B)$ to ${\cal N}(C,D)$ such that 
\[ \tilde\pi\rho_{A,B}(r)=\rho_{C,D}\pi(r) \]
for all $r\in{\cal M}(A,B)$. 

\begin{lem}\label{torsioniso}
The homomorphism $\tilde\pi$ induces an isomorphism 
between $T({\cal N}(A,B))$ and $T({\cal N}(C,D))$. 
\end{lem}
\begin{proof}
By using $\{u(c,a)\}_{c\in C,a\in A}$ and $\{v(d,b)\}_{d\in D,b\in B}$ 
such that 
\[ \sum_{c\in C}(e-\nu(c))u(c,a)=e-\mu(a) \]
and 
\[ \sum_{d\in D}(e-\nu(d))v(d,b)=e-\mu(b) \]
for all $a\in A$ and $b\in B$, 
we can define a homomorphism $\pi':{\cal M}_\Q(C,D)\to{\cal M}_\Q(A,B)$ 
in a similar fashion. 
It suffices to show that 
$\tilde\pi'\circ\tilde\pi$ induces the identity map 
on $T({\cal N}(A,B))$. 

Take $r\in\rho_{A,B}^{-1}(T({\cal N}(A,B)))$ arbitrarily and 
let $r'=\pi'(\pi(r))$. 
We would like to show $r'-r\in{\cal A}+{\cal B}$. 
There exist $p\in{\cal A}_\Q$ and $q\in{\cal B}_\Q$ 
such that $r=p+q$. 
We can find $\alpha:A\to\Q[G]$ and $\beta:B\to\Q[G]$ such that 
\[ p=\sum_{a\in A}\sum_{b\in B}
\alpha(a)(e-\mu(b))\otimes a\otimes b \]
and 
\[ q=\sum_{a\in A}\sum_{b\in B}
\beta(b)(e-\mu(a))\otimes a\otimes b. \]
Let $p'(a',b')\in\Q[G]$ be 
the $a'\otimes b'$-component of $\pi'(\pi(p))\in{\cal A}_\Q$, 
that is, 
\[ \pi'(\pi(p))=\sum_{a'\in A}\sum_{b'\in B}
p'(a',b')\otimes a'\otimes b'. \]
Then we have 
\begin{align*}
p'(a',b')
&=\sum_{a\in A}\sum_{b\in B}\sum_{c\in C}\sum_{d\in D}
\alpha(a)(e-\mu(b))s(a,c)t(b,d)u(c,a')v(d,b') \\
&=\sum_{a\in A}\sum_{c\in C}\alpha(a)(e-\mu(b'))s(a,c)u(c,a'), 
\end{align*}
and so 
\begin{align*}
p'(a',b')+\Z[G]
&=\sum_{a\in A}\sum_{c\in C}-\beta(b')(e-\mu(a))s(a,c)u(c,a')+\Z[G] \\
&=-\beta(b')(e-\mu(a'))+\Z[G] \\
&=\alpha(a')(e-\mu(b'))+\Z[G]. 
\end{align*}
Therefore we get 
\[ \pi'(\pi(p))-p\in{\cal M}(A,B). \]
By Lemma \ref{ABpure}, we obtain 
\[ \pi'(\pi(p))-p\in{\cal A}. \]
In the same way, 
\[ \pi'(\pi(q))-q\in{\cal B} \]
is obtained. 
Consequently, 
\[ \pi'(\pi(r))-r=(\pi'(\pi(p))-p)+(\pi'(\pi(q))-q)\in{\cal A}+{\cal B}. \]
Hence 
\[ \rho_{A,B}(r)=\rho_{A,B}(\pi'(\pi(r)))
=\tilde\pi'(\rho_{C,D}(\pi(r)))=\tilde\pi'(\tilde\pi(\rho_{A,B}(r))), \]
thereby completing the proof. 
\end{proof}

By the lemma above, 
we have noticed that 
$T({\cal N}(A,B))$ depends only on the finite abelian group $G$. 
Let us denote this finite abelian group $T({\cal N}(A,B))$ by $T_G$. 
\bigskip

Let us consider the case of $G=\Z_m$. 
Let $p\in\Z_m$ be a generator. 
Thanks to Lemma \ref{torsioniso}, we may assume that 
both $A$ and $B$ are singletons and $\mu(A)=\mu(B)=\{p\}$. 
Then we have ${\cal A}={\cal B}$. 
It follows from Lemma \ref{ABpure} that 
${\cal M}(A,B)/({\cal A}+{\cal B})={\cal M}(A,B)/{\cal A}$ is torsion free. 
Thus $T_G=T({\cal N}(A,B))$ is trivial. 
\bigskip

We now turn to the general case. 
Let $G=\Z_{m_1}\times\Z_{m_2}\times\dots\times\Z_{m_n}$ and 
let $p_i$ be the generator of the $i$-th component. 
Thus, 
\[ G=\{p_1^{k_1}p_2^{k_2}\dots p_n^{k_n}:
k_i=0,1,\dots,m_i-1\}. \]
We denote $e+p_i+p_i^2+\dots+p_i^{m_i-1}$ by $P_i$. 
Let 
\[ Q_i=\prod_{j\neq i}P_j \]
and let $N=\sum_{g\in G}g$. 
Notice that $N=P_1P_2\dots P_n$.  

\begin{lem}\label{solution}
Suppose that $\{r_i\}_{i=1}^n\subset\Q[G]$ satisfies 
$r_iP_i\in\Z[G]$ and 
\[ r_i(e-p_j)-r_j(e-p_i)\in\Z[G] \]
for all $i,j=1,2,\dots,n$. 
Then we can find $x\in\Q[G]$ such that 
\[ r_i-x(e-p_i)\in\Z[G] \]
for all $i=1,2,\dots,n$. 
\end{lem}
\begin{proof}
The proof goes by induction on $n$. 
Let us consider the case of $n=1$. 
We can write 
\[ r_1=\sum_{j=1}^{m_1-1}s_jp_1^j. \]
Define 
\[ x=\sum_{j=1}^{m_1-1}\sum_{k=1}^js_kp_1^j. \]
Then $r_1P_1\in\Z[G]$ ensures 
\[ r_1-x(e-p_1)\in\Z[G]. \]

Suppose that the assertion has been proved in the case of $n-1$. 
Let $H$ be the subgroup generated by $p_2,p_3,\dots,p_n$ and 
let $G_1$ be the subgroup generated by $p_1$. 
We regard $\Q[H]$ and $\Q[G_1]$ as subrings of $\Q[G]$. 
For each $i=2,3,\dots,n$, 
we can find $r_{i,0},r_{i,1},\dots,r_{i,m_1-1}\in\Q[H]$ such that 
\[ r_i=\sum_{k=0}^{m_1-1}r_{i,k}p_1^k. \]
It is easy to see $r_{i,k}P_i\in\Z[H]$ and 
\[ r_{i,k}(e-p_j)-r_{j,k}(e-p_i)\in\Z[H] \]
for all $i,j=2,3,\dots,n$ and $k=0,1,\dots,m_1-1$. 
For each $k$, by the induction hypothesis, 
there exists $y_k\in\Q[H]$ such that 
\[ r_{i,k}-y_k(e-p_i)\in\Z[H] \]
for all $i=2,3,\dots,n$. 
Put 
\[ y=\sum_{k=0}^{m_1-1}y_kp_1^k. \]
Then it is routine to check 
\[ r_i-y(e-p_i)\in\Z[G] \]
for all $i=2,3,\dots,n$. 
Let us consider $r_1-y(e-p_1)$. 
For any $i=2,3,\dots,n$, we have  
\[ (r_1-y(e-p_1))(e-p_i)\in\Z[G], \]
which means that there exists $s\in\Q[G_1]$ such that 
\[ r_1-y(e-p_1)-sQ_1\in\Z[G]. \]
Moreover, $r_1P_1\in\Z[G]$ implies $sP_1\in\Z[G_1]$. 
By the induction hypothesis for $n=1$, 
we can find $z\in\Q[G_1]$ such that 
\[ s-z(e-p_1)\in\Z[G_1]. \]
We define $x\in\Q[G]$ by 
\[ x=y+zQ_1. \]
It is clear that $r_i-x(e-p_i)$ belongs to $\Z[G]$ 
for all $i=2,3,\dots,n$. 
Besides, we can see 
\[ r_1-x(e-p_1)=r_1-y(e-p_1)-z(e-p_1)Q_1\in\Z[G], \]
which completes the proof. 
\end{proof}

\begin{lem}\label{adjust}
Let $r_1,r_2,\dots,r_n\in\Q[G]$. 
Suppose 
\[ r_i(e-p_j)-r_j(e-p_i)\in\Z[G] \]
for all $i,j=1,2,\dots,n$. 
Then there exist $u_1,u_2,\dots,u_n\in\Z[G]$ such that 
\[ (r_i-u_i)(e-p_j)=(r_j-u_j)(e-p_i) \]
for all $i,j=1,2,\dots,n$. 
\end{lem}
\begin{proof}
Let us consider $r_iP_i$. 
For any $j=1,2,\dots,n$, we have 
\[ r_iP_i(e-p_j)+\Z[G]=r_jP_i(e-p_i)+\Z[G]=\Z[G]. \]
It follows that there exists $t_i\in\Q$ such that 
\[ r_iP_i-t_iN\in\Z[G]. \]
Define $s_i\in\Q[G]$ by 
\[ s_i=r_i-t_iQ_i. \]
It is easy to see $s_iP_i\in\Z[G]$ and 
\[ s_i(e-p_j)-s_j(e-p_i)\in\Z[G] \]
for all $i,j=1,2,\dots,n$. 
By Lemma \ref{solution}, there exists $x\in\Q[G]$ such that 
\[ s_i-x(e-p_i)\in\Z[G]. \]
It is not hard to see that $u_i=s_i-x(e-p_i)$ does the work. 
\end{proof}

We would like to compute $T_G\cong T({\cal N}(A,B))$. 
Thanks to Lemma \ref{torsioniso}, 
we may assume that $A=B=\{p_1,p_2,\dots,p_n\}$ and that 
the map $\mu$ is the inclusion to $G$. 
To simplify the notation, we omit $\mu$. 

Suppose that 
\[ r=\sum_{i,j=1}^nr(p_i,p_j)\otimes p_i\otimes p_j
\in{\cal M}_\Q(A,B) \]
belongs to ${\cal A}_\Q+{\cal B}_\Q$. 
There exist $\alpha,\beta:\{p_1,p_2,\dots,p_n\}\to\Q[G]$ such that 
\[ r(a,b)=\alpha(a)(e-b)+\beta(b)(e-a) \]
for all $a,b\in\{p_1,p_2,\dots,p_n\}$. 
For every $i,j=1,2,\dots,n$, 
let us consider the element 
\[ s_{i,j}=(\alpha(p_i)+\beta(p_i))(e-p_j). \]
In order to show that 
$s_{i,j}$ does not depend on the choice of $\alpha$ and $\beta$, 
suppose that there exist $\alpha',\beta':\{p_1,p_2,\dots,p_n\}\to\Q[G]$ 
such that 
\[ r(a,b)=\alpha'(a)(e-b)+\beta'(b)(e-a) \]
for all $a,b\in\{p_1,p_2,\dots,p_n\}$. 
As before, we put 
\[ s'_{i,j}=(\alpha'(p_i)+\beta'(p_i))(e-p_j). \]
For any $k=1,2,\dots,n$, we have 
\begin{align*}
& s_{i,j}(e-p_k) \\
&=\alpha(p_i)(e-p_j)(e-p_k)+r(p_k,p_i)(e-p_j)-\alpha(p_k)(e-p_j)(e-p_i) \\
&=\alpha(p_i)(e-p_j)(e-p_k)+r(p_k,p_i)(e-p_j)
-r(p_k,p_j)(e-p_i)+\beta(p_j)(e-p_k)(e-p_i) \\
&=r(p_k,p_i)(e-p_j)-r(p_k,p_j)(e-p_i)+r(p_i,p_j)(e-p_k), 
\end{align*}
and so 
\[ s_{i,j}(e-p_k)=s'_{i,j}(e-p_k). \]
Combining this with $(s_{i,j}-s'_{i,j})P_j=0$, we obtain $s_{i,j}=s'_{i,j}$. 
Thus, $s_{i,j}$ does not depend on the choice of $\alpha$ and $\beta$. 
We write $\kappa_{i,j}(r)=s_{i,j}$. 
Clearly $\kappa_{i,j}$ is a homomorphism 
from ${\cal A}_\Q+{\cal B}_\Q$ to $\Q[G]$. 

Next, suppose that 
\[ r=\sum_{i,j=1}^nr(p_i,p_j)\otimes p_i\otimes p_j
\in{\cal M}(A,B) \]
belongs to $\rho_{A,B}^{-1}(T({\cal N}(A,B)))$. 
As mentioned before, 
this means that $r$ is in ${\cal A}_\Q+{\cal B}_\Q$. 

For distinct $i,j\in\{1,2,\dots,n\}$, 
let us consider $\kappa_{i,j}(r)$. 
For any $k=1,2,\dots,n$, we have 
\begin{align*}
& \kappa_{i,j}(r)(e-p_k)+\Z[G] \\
&=(\alpha(p_i)+\beta(p_i))(e-p_j)(e-p_k)+\Z[G] \\
&=\alpha(p_i)(e-p_j)(e-p_k)-\alpha(p_k)(e-p_j)(e-p_i)+\Z[G] \\
&=\alpha(p_i)(e-p_j)(e-p_k)+\beta(p_j)(e-p_k)(e-p_i)+\Z[G] \\
&=r(p_i,p_j)(e-p_k)+\Z[G]=\Z[G]. 
\end{align*}
Therefore there exists $t_{i,j}\in\Q$ such that 
\[ \kappa_{i,j}(r)-t_{i,j}N\in\Z[G]. \]
Because $t_{i,j}$ is uniquely determined up to $\Z$, 
we can define a map $\pi_{i,j}$ 
sending $r$ to $e^{2\pi\sqrt{-1}t_{i,j}}\in\T$. 
It is not hard to see that $\pi_{i,j}$ is a homomorphism 
from $\rho_{A,B}^{-1}(T({\cal N}(A,B)))$ to $\T$. 
We remark that $\pi_{j,i}(r)=\pi_{i,j}(r)^{-1}$. 

By the definition of $\kappa_{i,j}(r)$, 
we have $\kappa_{i,j}(r)P_j=0$. 
Then we have $t_{i,j}NP_j=t_{i,j}m_jN\in\Z[G]$, 
which implies $t_{i,j}m_j\in\Z$. 
It follows from 
\[ \kappa_{i,j}(r)
=r(p_i,p_j)+r(p_j,p_i)-(\alpha(p_j)+\beta(p_j))(e-p_i), \]
that $\kappa_{i,j}(r)P_i$ belongs to $\Z[G]$. 
Hence $t_{i,j}NP_i=t_{i,j}m_iN\in\Z[G]$, 
which means $t_{i,j}m_i\in\Z$. 
Consequently, we obtain 
\[ \pi_{i,j}(r)=e^{2\pi\sqrt{-1}t_{i,j}}
\in\{\omega\in\T:\omega^{d(i,j)}=1\}, \]
where $d(i,j)$ is the greatest common divisor of $m_i$ and $m_j$. 
Thus, we regard $\pi_{i,j}$ as a homomorphism 
from $\rho_{A,B}^{-1}(T({\cal N}(A,B)))$ to $\Z_{d(i,j)}$. 

We define a homomorphism 
\[ \pi:\rho_{A,B}^{-1}(T({\cal N}(A,B)))\to\prod_{i<j}\Z_{d(i,j)} \]
by $\pi(r)=(\pi_{i,j}(r))_{i,j}$. 

\begin{lem}
The homomorphism $\pi$ is surjective. 
\end{lem}
\begin{proof}
Take distinct $k,l\in\{1,2,\dots,n\}$ arbitrarily. 
Put $d=d(k,l)$. 
It suffices to show there exists $r\in\rho_{A,B}^{-1}(T({\cal N}(A,B)))$ 
such that 
\[ \pi_{i,j}(r)
=\begin{cases}e^{2\pi\sqrt{-1}/d} & \text{ if }\{i,j\}=\{k,l\} \\
1 & \text{ if }\{i,j\}\neq\{k,l\}. \end{cases} \]
Let 
\[ r=\left(\frac{m_k}{d}Q_k-\frac{m_l}{d}Q_l\right)
\otimes p_k\otimes p_l\in\Z[G]. \]
We define $\alpha,\beta:\{p_1,p_2,\dots,p_n\}\to\Q[G]$ by 
\[ \alpha(p_i)
=\begin{cases}\displaystyle\sum_{j=0}^{m_l-1}\frac{j}{d}p_l^jQ_l 
& \text{ if }i=k \\
0 & \text{ if }i\neq k, \end{cases} \]
and 
\[ \beta(p_i)
=\begin{cases}\displaystyle\sum_{j=0}^{m_k-1}-\frac{j}{d}p_k^jQ_k 
& \text{ if }i=l \\
0 & \text{ if }i\neq l. \end{cases} \]
We can easily check 
\[ r=\sum_{i,j=1}^n
(\alpha(p_i)(e-p_j)+\beta(p_j)(e-p_i))\otimes p_i\otimes p_j. \]
In addition, we get 
\[ \kappa_{k,l}(r)
=(\alpha(p_k)+\beta(p_k))(e-p_l)
=\frac{1}{d}N-\frac{m_l}{d}Q_l. \]
It follows that 
\[ \kappa_{k,l}(r)-\frac{1}{d}N\in\Z[G]. \]
Therefore we obtain $\pi_{k,l}(r)=e^{2\pi\sqrt{-1}/d}$. 
If $\{i,j\}\neq\{k,l\}$, then $\kappa_{i,j}(r)=0$, 
and so $\pi_{i,j}(r)=1$. 
\end{proof}

\begin{lem}
The kernel of $\pi$ is equal to ${\cal A}+{\cal B}$. 
\end{lem}
\begin{proof}
Evidently ${\cal A}+{\cal B}$ is contained in $\ker\pi$. 
Suppose that 
\[ r=\sum_{i,j=1}^nr(p_i,p_j)\otimes p_i\otimes p_j\in{\cal M}(A,B) \]
belongs to $\ker\pi$. 
We would like to show $r\in{\cal A}+{\cal B}$. 
There exist $\alpha,\beta:\{p_1,p_2,\dots,p_n\}\to\Q[G]$ such that 
\[ r(a,b)=\alpha(a)(e-b)+\beta(b)(e-a) \]
for all $a,b\in\{p_1,p_2,\dots,p_n\}$. 
By assumption, we have 
\[ (\alpha(p_i)+\beta(p_i))(e-p_j)\in\Z[G] \]
for all $i,j=1,2,\dots,n$. 
This implies that there exists $k_i\in\Q$ such that 
\[ \alpha(p_i)+\beta(p_i)\in k_iN+\Z[G]. \]
By replacing $\alpha(p_i)$ with $\alpha(p_i)-k_iN$, 
we may assume that $\alpha(p_i)+\beta(p_i)\in\Z[G]$ 
for all $i=1,2,\dots,n$. 
Then we obtain 
\[ \alpha(p_i)(e-p_j)-\alpha(p_j)(e-p_i)
=\alpha(p_i)(e-p_j)+\beta(p_i)(e-p_j)-r(p_j,p_i)\in\Z[G]. \]
By applying Lemma \ref{adjust} to 
$\alpha(p_1),\alpha(p_2),\dots,\alpha(p_n)$, 
we obtain $u_1,u_2,\dots,u_n\in\Z[G]$ such that 
\[ (\alpha(p_i)-u_i)(e-p_j)=(\alpha(p_j)-u_j)(e-p_i) \]
for all $i,j=1,2,\dots,n$. 
We define $\tilde\alpha,\tilde\beta:\{p_1,p_2,\dots,p_n\}\to\Z[G]$ by 
\[ \tilde\alpha(p_i)=u_i \]
and 
\[ \tilde\beta(p_i)=\alpha(p_i)+\beta(p_i)-u_i. \]
It is easy to check 
\[ r(p_i,p_j)=\tilde\alpha(p_i)(e-p_j)+\tilde\beta(p_j)(e-p_i) \]
for all $i,j=1,2,\dots,n$. 
Therefore $r$ belongs to ${\cal A}+{\cal B}$. 
\end{proof}

We can summarize the results obtained above as follows. 

\begin{thm}\label{T_G=wedge}
Let $G=\Z_{m_1}\times\Z_{m_2}\times\dots\times\Z_{m_n}$ and 
let $d(i,j)$ be the greatest common divisor of $m_i$ and $m_j$. 
Then $T_G$ is isomorphic to the direct product of 
all $\Z_{d(i,j)}$'s for $1\leq i<j\leq n$. 
In other words $T_G$ is isomorphic to the wedge product $G\wedge G$. 
\end{thm}

\section{Skew product extensions of product systems}

Let $\xp$ be a Cantor minimal $\Z$-system. 
The group of coinvariants of $\xp$ 
\[ \Coinv\xp=C(X,\Z)/\{f-f\circ\phi^{-1}:f\in C(X,\Z)\}. \]
was denoted by $K^0\xp$ in \cite{HPS}. 
At first, we must recall that 
$\Coinv\xp$ can be written as an inductive limit of 
free abelian groups of finite rank. 
The reader may refer to \cite{HPS} for the details. 

A family of non-empty clopen subsets 
\[ {\cal P}=\{X(v,k):v\in V,1\leq k\leq h(v)\} \]
indexed by a finite set $V$ and natural numbers $k=1,2,\dots,h(v)$ 
is called a Kakutani-Rohlin partition of $\xp$, 
if the following conditions are satisfied: 
\begin{itemize} 
\item ${\cal P}$ is a partition of $X$. 
\item For all $v\in V$ and $k=1,2,\dots,h(v)-1$, we have 
$\phi(X(v,k))=X(v,k+1)$. 
\end{itemize}
We define the roof set $R({\cal P})$ and the base set $B({\cal P})$ by 
\[ R({\cal P})=\bigcup_{v\in V}X(v,h(v)) \]
and 
\[ B({\cal P})=\bigcup_{v\in V}X(v,1). \]
Notice that $\phi(R({\cal P}))=B({\cal P})$. 
For each $v\in V$, the family of clopen sets 
$X(v,1),X(v,2),\dots,X(v,h(v))$ is called a tower, 
and $h(v)$ is called the height of the tower. 
We may identify the label $v$ with the corresponding tower. 

Let $\{{\cal P}_n\}_{n\in\N}$ be a sequence of Kakutani-Rohlin partitions. 
We denote the set of towers in ${\cal P}_n$ by $V_n$ 
and clopen sets belonging to ${\cal P}_n$ 
by $X(n,v,k)$ ($v\in V_n,k=1,2,\dots,h(v)$). 
We say that 
$\{{\cal P}_n\}_{n\in\N}$ gives a Bratteli-Vershik model for $\xp$, 
if the following are satisfied: 
\begin{itemize}
\item The roof sets $R({\cal P}_n)=\bigcup_{v\in V_n}X(n,v,h(v))$ 
form a decreasing sequence of clopen sets, 
which shrinks to a single point. 
\item ${\cal P}_{n+1}$ is finer than ${\cal P}_n$ for all $n\in\N$ 
as partitions, and $\bigcup_n{\cal P}_n$ generates the topology of $X$. 
\end{itemize}
Let $\Z^{V_n}$ be the free abelian group which has basis $\{v:v\in V_n\}$. 
For any $v\in V_n$ and $v'\in V_{n+1}$, we put 
\[ A_n(v,v')=\#\{1\leq k\leq h(v'):X(n+1,v',k)\subset X(n,v,1)\}. \]
We can define a homomorphism $\alpha_n:\Z^{V_n}\to\Z^{V_{n+1}}$ by 
\[ \alpha_n(v)=\sum_{v'\in V_{n+1}}A_n(v,v')v'\in\Z^{V_{n+1}} \]
for all $v\in V_n$. 

\begin{thm}
The group of coinvariants $\Coinv\xp$ is 
isomorphic to the inductive limit of $\Z^{V_n}$ via $\alpha_n$, that is, 
$\Coinv\xp\cong\lim_{n\to\infty}(\Z^{V_n},\alpha_n)$. 
\end{thm}
\bigskip

We would like to consider a skew product extension of $\xp$ 
associated with a finite group valued cocycle. 
Let $G$ be a finite group and let $\xi:X\to G$ be a continuous map. 
We call $\xi$ a cocycle taking its values in $G$. 
Define a homeomorphism $\phi\times\xi$ on $X\times G$ by 
\[ \phi\times\xi(x,g)=(\phi(x),g\xi(x)) \]
for all $x\in X$ and $g\in G$. 
The cocycle $\xi$ is said to be non-degenerate, 
if $\phi\times\xi$ is a minimal homeomorphism. 

Given a cocycle $\xi$, we may assume that 
$\xi$ is constant on every clopen set $U\in{\cal P}_1$. 
Since ${\cal P}_n$ is finer than ${\cal P}_1$, 
this means that $\xi$ is constant 
on every clopen set $U\in{\cal P}_n$ for all $n\in\N$. 
For any $n\in\N$ and $v\in V_n$, we define 
\[ \tilde\xi(v)=\xi(x)\xi(\phi(x))\xi(\phi^2(x))
\dots\xi(\phi^{h(v)-1}(x)), \]
where $x$ is an arbitrary point in $X(n,v,1)$. 
It is easy to see that $\tilde\xi(v)$ does not depend on the choice of $x$, 
because $\xi$ is constant on each $X(n,v,k)$. 
The next lemma is quoted from \cite{M}. 

\begin{lem}[{\cite[Lemma 2.3]{M}}]\label{nondege}
The cocycle $\xi$ is non-degenerate if and only if 
$\{\tilde\xi(v):v\in V_n\}$ generates $G$ for every $n\in\N$. 
\end{lem}
\bigskip

From now on, we assume that $G$ is a finite abelian group. 
We still use the multiplicative notation. 
Let $\xp$ and $\yp$ be two Cantor minimal $\Z$-systems and 
let $\xi:X\to G$ and $\eta:Y\to G$ be two non-degenerate cocycles. 

Put $Z=X\times Y\times G$. 
We define $\tilde\phi,\tilde\psi\in\Homeo(Z)$ by 
\[ \tilde\phi(x,y,g)=(\phi(x),y,g\xi(x)) \]
and 
\[ \tilde\psi(x,y,g)=(x,\psi(y),g\eta(y)). \]
Evidently $\tilde\phi$ commutes with $\tilde\psi$, 
and so they form a $\Z^2$-action $\omega$ on $Z$. 
It is easy to see $(Z,\omega)$ is a minimal dynamical system. 
The following is the main theorem of this paper. 

\begin{thm}\label{main}
The torsion part of the group of coinvariants $\Coinv(Z,\omega)$ is 
isomorphic to $T_G$. 
\end{thm}

For example, if $G$ is the direct product of $n$ copies of $\Z_m$, 
then the torsion subgroup of $\Coinv(Z,\omega)$ is isomorphic to 
the direct product of $n(n-1)/2$ copies of $\Z_m$. 

We need a series of lemmas to prove this theorem. 

The group of coinvariants $\Coinv(Z,\omega)$ is the quotient 
of $C(Z,\Z)$ by 
\[ D=\{f_1-f_1\circ\tilde\phi^{-1}+f_2-f_2\circ\tilde\psi^{-1}:
f_1,f_2\in C(Z,\Z)\}. \]
Fix $x_0\in X$ and $y_0\in Y$. 
Let $D'\subset D$ be the subgroup consisting of 
all elements $f_1-f_1\circ\tilde\phi^{-1}+f_2-f_2\circ\tilde\psi^{-1}$ 
with $f_1(x_0,y,g)=f_2(x,y_0,g)=0$ for all $(x,y,g)\in Z$. 
At first, we would like to consider 
the quotient of $C(Z,\Z)$ by $D'$. 

We take Bratteli-Vershik models 
$\{{\cal P}_n\}_{n\in\N}$ and $\{{\cal Q}_n\}_{n\in\N}$ 
for $\xp$ and $\yp$, respectively. 
Let $V_n$ and $W_n$ denote the towers 
in ${\cal P}_n$ and ${\cal Q}_n$, respectively. 
We write 
\[ {\cal P}_n=\{X(n,v,k):v\in V_n,1\leq k\leq h(v)\} \]
and 
\[ {\cal Q}_n=\{Y(n,w,l):w\in W_n,1\leq l\leq h(w)\} \]
for every $n\in\N$. 
We may assume $\bigcap_nR({\cal P}_n)=\{x_0\}$ 
and $\bigcap_nR({\cal Q}_n)=\{y_0\}$. 
Moreover, we may also assume that 
$\xi$ (resp. $\eta$) is constant 
on every clopen set in ${\cal P}_1$ (resp. ${\cal Q}_1$). 
Then $\tilde\xi(v)$ and $\tilde\eta(w)$ are well-defined 
for every $v\in V_n$ and $w\in W_n$. 
Besides, we define $\tilde\xi(v,k)$ by $\tilde\xi(v,1)=e$ and by 
\[ \tilde\xi(v,k)=\xi(x)\xi(\phi(x))\xi(\phi^2(x))
\dots\xi(\phi^{k-2}(x)) \]
for $k=2,3,\dots,h(v)$, 
where $x$ is a point in $X(n,v,1)$. 
For every $w\in W_n$ and $l=1,2,\dots,h(w)$, 
$\tilde\eta(w,l)$ can be defined in a similar fashion. 

Let $\chi(v,w,k,l,g)\in C(Z,\Z)$ be the characteristic function 
on $X(n,v,k)\times Y(n,w,l)\times\{g\}$ and 
let $C_n\subset C(Z,\Z)$ be the subgroup generated by 
\[ \{\chi(v,w,k,l,g):
v\in V_n,w\in W_n,1\leq k\leq h(v),1\leq l\leq h(w),g\in G\}. \]
Evidently $C_n$ is a finitely generated free abelian group and 
$C(Z,\Z)$ is equal to $\bigcup_n C_n$. 
Let $D'_n\subset D'$ be the subgroup consisting of 
all elements $f_1-f_1\circ\tilde\phi^{-1}+f_2-f_2\circ\tilde\psi^{-1}$ 
with $f_1,f_2\in C_n$ and 
\[ f_1|_{R({\cal P}_n)\times Y\times G}
=f_2|_{X\times R({\cal Q}_n)\times G}=0. \]
Since $\xi$ (resp. $\eta$) is constant on each clopen set 
belonging to ${\cal P}_n$ (resp. ${\cal Q}_n$), 
we can see that $D'_n$ is contained in $C_n$. 
Clearly we have $D'_n\subset D'_{n+1}$ 
and $D'=\bigcup_nD'_n$. 
It follows that $C(Z,\Z)/D'$ is isomorphic to 
the inductive limit of $C_n/D'_n$'s. 

It is easily verified that 
\[ \chi(v,w,1,1,g)-\chi(v,w,k,l,g\tilde\xi(v,k)\tilde\eta(w,l)) \]
belongs to $D'_n$ and $D'_n$ is spanned by these elements. 
We put ${\cal M}(V_n,W_n)=\Z[G]\otimes\Z^{V_n}\otimes\Z^{W_n}$ and 
denote each element $r$ of ${\cal M}(V_n,W_n)$ by 
\[ r=\sum_{v\in V_n}\sum_{w\in W_n}
r(v,w)\otimes v\otimes w, \]
where $r(v,w)$ is an element of the group ring $\Z[G]$. 
Define $\tau_n:C_n\to{\cal M}(V_n,W_n)$ by 
\[ \tau_n(\chi(v,w,k,l,g))
=g\tilde\xi(v,k)^{-1}\tilde\eta(w,l)^{-1}\otimes v\otimes w. \]
Then the kernel of $\tau_n$ agrees with $D'_n$. 
Let $\pi_n:{\cal M}(V_n,W_n)\to{\cal M}(V_{n+1},W_{n+1})$ be 
the connecting homomorphism induced by the inclusion $C_n\subset C_{n+1}$. 
Thus, we get the following commutative diagram: 
\[ \begin{CD}
0@>>>D'_n@>>>C_n@>\tau_n>>{\cal M}(V_n,W_n)@>>>0 \\
@.@VVV@VVV@VV\pi_nV@. \\
0@>>>D'_{n+1}@>>>C_{n+1}@>\tau_{n+1}>>{\cal M}(V_{n+1},W_{n+1})@>>>0, 
\end{CD} \]
where the vertical arrows from $D'_n$ to $D'_{n+1}$ and 
from $C_n$ to $C_{n+1}$ are inclusions. 
As a direct consequence, we have the following. 

\begin{lem}\label{AFpart}
In this setting, $C(Z,\Z)/D'$ is isomorphic to 
the inductive limit of 
${\cal M}(V_n,W_n)=\Z[G]\otimes\Z^{V_n}\otimes\Z^{W_n}$. 
In particular, $C(Z,\Z)/D'$ is torsion free. 
\end{lem}

We would like to see the connecting homomorphism $\pi_n$ 
from ${\cal M}(V_n,W_n)$ to ${\cal M}(V_{n+1},W_{n+1})$. 
For $v\in V_n$ and $v'\in V_{n+1}$, we set 
\[ E_n(v,v')=\{1\leq k\leq h(v'):X(n+1,v',k)\subset X(n,v,1)\} \]
and 
\[ s(v,v')=\sum_{k\in E_n(v,v')}\tilde\xi(v',k)^{-1}\in\Z[G]. \]
In the same way, $t(w,w')\in\Z[G]$ can be defined 
for $w\in W_n$ and $w'\in W_{n+1}$. 
Since 
\[ \chi(v',w',1,1,g\tilde\xi(v',k)^{-1}\tilde\eta(w',l)^{-1})
-\chi(v',w',k,l,g) \]
belongs to $D'_{n+1}$, we obtain 
\[ \pi_n(r\otimes v\otimes w)
=\sum_{v'\in V_{n+1}}\sum_{w'\in W_{n+1}}
rs(v,v')t(w,w')\otimes v'\otimes w' \]
for all $r\in\Z[G]$. 

\begin{lem}
\begin{enumerate}\label{connect}
\item For every $v'\in V_{n+1}$, we have 
\[ \sum_{v\in V_n}(e-\tilde\xi(v)^{-1})s(v,v')=e-\tilde\xi(v')^{-1}. \]
\item For every $w'\in W_{n+1}$, we have 
\[ \sum_{w\in W_n}(e-\tilde\eta(w)^{-1})t(w,w')=e-\tilde\eta(w')^{-1}. \]
\end{enumerate}
\end{lem}
\begin{proof}
It suffices to show (1). 
Put 
\[ E=\bigcup_{v\in V_n}E_n(v,v'). \]
It is clear that 
\[ \sum_{v\in V_n}s(v,v')=\sum_{k\in E}\tilde\xi(v',k)^{-1}. \]
Note that $E$ contains one, because 
\[ X(n+1,v',1)\subset B({\cal P}_{n+1})
\subset B({\cal P}_n)=\bigcup_{v\in V_n}X(n,v,1). \]
Let $1=k_1<k_2<\dots<k_L$ be the arranged list of elements in $E$. 
If $k_i$ belongs to $E_n(v,v')$ and $i\neq L$, 
then we can see that 
\[ \tilde\xi(v)^{-1}\tilde\xi(v',k_i)^{-1}=\tilde\xi(v',k_{i+1})^{-1}. \]
If $k_L$ belongs to $E_n(v,v')$, then 
we get 
\[ \tilde\xi(v)^{-1}\tilde\xi(v',k_L)^{-1}=\tilde\xi(v')^{-1}. \]
It follows that 
\[ \sum_{v\in V_n}\tilde\xi(v)^{-1}s(v,v')
=\sum_{i=2}^L\tilde\xi(v',k_i)^{-1}+\tilde\xi(v')^{-1}. \]
Hence we obtain 
\[ \sum_{v\in V_n}(e-\tilde\xi(v)^{-1})s(v,v')
=\tilde\xi(v',1)-\tilde\xi(v')^{-1}=e-\tilde\xi(v')^{-1}. \]
\end{proof}

For $v\in V_n$ and $g\in G$, 
we define $\partial_Y(v,g)\in C_n$ by 
\[ \partial_Y(v,g)=\sum_{w\in W_n}\chi(v,w,1,1,g)
-\sum_{w\in W_n}\chi(v,w,1,1,g)\circ\tilde\psi. \]
Similarly, for $w\in W_n$ and $g\in G$, 
we define $\partial_X(w,g)\in C_n$ by 
\[ \partial_X(w,g)=\sum_{v\in V_n}\chi(v,w,1,1,g)
-\sum_{v\in V_n}\chi(v,w,1,1,g)\circ\tilde\phi. \]
Let $D_n$ be the subgroup of $C_n$ generated by 
$D'_n$ and 
\[ \{\partial_X(w,g):w\in W_n,g\in G\}
\cup\{\partial_Y(v,g):v\in V_n,g\in G\}. \]

\begin{lem}
\begin{enumerate}
\item For any $g\in G$, $v\in V_n$ and $k=1,2,\dots,h(v)$, 
\[ \sum_{w\in W_n}\chi(v,w,k,1,g)
-\sum_{w\in W_n}\chi(v,w,k,1,g)\circ\tilde\psi \]
is contained in $D_n$. 
\item For any $g\in G$, $w\in W_n$ and $l=1,2,\dots,h(w)$, 
\[ \sum_{v\in V_n}\chi(v,w,1,l,g)
-\sum_{v\in V_n}\chi(v,w,1,l,g)\circ\tilde\psi \]
is contained in $D_n$. 
\end{enumerate}
\end{lem}
\begin{proof}
It suffices to show (1). 
At first, we have 
\[ \chi(v,w,k,1,g)-\chi(v,w,1,1,g\tilde\xi(v,k)^{-1})\in D'_n \]
for every $w\in W_n$. It follows that 
\[ \sum_{w\in W_n}\chi(v,w,k,1,g)
-\sum_{w\in W_n}\chi(v,w,1,1,g\tilde\xi(v,k)^{-1})\in D'_n. \]
Put $\eta_w=\tilde\eta(w)^{-1}\tilde\eta(w,h(w))$ for every $w\in W_n$. 
It is not hard to see 
\[ \sum_{w\in W_n}\chi(v,w,k,1,g)\circ\tilde\psi
=\sum_{w\in W_n}\chi(v,w,k,h(w),g\eta_w) \]
for every $g\in G$. 
By combining this with 
\[ \chi(v,w,k,h(w),g\eta_w)
-\chi(v,w,1,h(w),g\eta_w\tilde\xi(v,k)^{-1})\in D'_n, \]
we obtain 
\[ \sum_{w\in W_n}\chi(v,w,k,1,g)\circ\tilde\psi
-\sum_{w\in W_n}\chi(v,w,1,1,g\tilde\xi(v,k)^{-1})\circ\tilde\psi\in D'_n. \]
Since 
\[ \sum_{w\in W_n}\chi(v,w,1,1,g\tilde\xi(v,k)^{-1})
-\sum_{w\in W_n}\chi(v,w,1,1,g\tilde\xi(v,k)^{-1})\circ\tilde\psi
=\partial_Y(v,g\tilde\xi(v,k)^{-1}), \]
we can conclude 
\[ \sum_{w\in W_n}\chi(v,w,k,1,g)
-\sum_{w\in W_n}\chi(v,w,k,1,g)\circ\tilde\psi\in D_n. \]
\end{proof}

\begin{lem}
In the setting above, 
we have $D_n\subset D_{n+1}$ and $D=\bigcup_nD_n$. 
\end{lem}
\begin{proof}
For any $v\in V_n$, 
$X(n,v,1)$ is a union of some clopen sets in ${\cal P}_{n+1}$, 
because ${\cal P}_{n+1}$ is finer than ${\cal P}_n$. 
It follows from (1) of the lemma above that 
$\partial_Y(v,g)$ is contained in $D_{n+1}$ for every $g\in G$. 
In the same way, 
$\partial_X(w,g)$ is contained in $D_{n+1}$ 
for every $w\in W_n$ and $g\in G$. 
Therefore we get $D_n\subset D_{n+1}$. 

Let us prove $D=\bigcup_nD_n$. 
It suffices to show that 
$1_O-1_O\circ\tilde\phi$ and $1_O-1_O\circ\tilde\psi$ belong 
to $\bigcup_nD_n$ for any clopen set $O\subset Z$. 
We may assume that there exist 
$g\in G$ and clopen sets $U\subset X$, $V\subset Y$ 
such that $O=U\times V\times\{g\}$. 

If $V$ does not contain $\psi(y_0)$, 
it is clear that $1_O-1_O\circ\tilde\psi$ belongs to $D'=\bigcup_nD'_n$. 
Then we have $1_O-1_O\circ\tilde\psi\in\bigcup_nD_n$, 
because $D'_n\subset D_n$. 
Suppose that $V$ contains $\psi(y_0)$. 
Since $\bigcup_n{\cal P}_n$ generates the topology of $X$ and 
$\bigcap_nB({\cal Q}_n)=\psi(y_0)$, 
there exists $n\in\N$ such that $1_U\in C_n$ and $B({\cal Q}_n)\subset V$. 
Put $O'=U\times B({\cal Q}_n)\times\{g\}$ and $O''=O\setminus O'$. 
We have already seen that 
$1_{O''}-1_{O''}\circ\tilde\psi$ is contained in $\bigcup_nD_n$, 
and so it suffices to show $1_{O'}-1_{O'}\circ\tilde\psi\in\bigcup_nD_n$. 
By $1_U\in C_n$, 
the clopen set $U$ is a union of some clopen sets in ${\cal P}_n$. 
For any $X(n,v,k)\in{\cal P}_n$, 
the characteristic function on $X(n,v,k)\times B({\cal Q}_n)\times\{g\}$ 
equals $\sum_{w\in W_n}\chi(v,w,k,1,g)$. 
Hence, from (1) of the lemma above, we can conclude 
$1_{O'}-1_{O'}\circ\tilde\psi\in\bigcup_nD_n$. 

In the same way we can prove $1_O-1_O\circ\tilde\phi\in\bigcup_nD_n$, 
and so we obtain $D=\bigcup_nD_n$. 
\end{proof}

For any $v\in V_n$ and $g\in G$, it is easy to see 
\[ \tau_n(\partial_Y(v,g))
=\sum_{w\in W_n}g(e-\tilde\eta(w)^{-1})\otimes v\otimes w. \]
Similarly, for any $w\in W_n$ and $g\in G$, 
\[ \tau_n(\partial_X(w,g))
=\sum_{v\in V_n}g(e-\tilde\xi(v)^{-1})\otimes v\otimes w. \]
Let ${\cal V}_n\subset\Z[G]\otimes\Z^{V_n}\otimes\Z^{W_n}$ be 
the subgroup consisting of 
\[ \sum_{v\in V_n}\sum_{w\in W_n}
\alpha(v)(e-\tilde\eta(w)^{-1})\otimes v\otimes w \]
with some map $\alpha:V_n\to\Z[G]$. 
In the same way, 
let ${\cal W}_n\subset\Z[G]\otimes\Z^{V_n}\otimes\Z^{W_n}$ be 
the subgroup consisting of 
\[ \sum_{v\in V_n}\sum_{w\in W_n}
\beta(w)(e-\tilde\xi(v)^{-1})\otimes v\otimes w \]
with some map $\beta:W_n\to\Z[G]$. 
Then it is evident that 
$\tau_n(D_n)$ coincides with ${\cal V}_n+{\cal W}_n$, 
that is, $C_n/D_n$ is isomorphic to 
${\cal M}(V_n,W_n)/({\cal V}_n+{\cal W}_n)$. 
By defining the map $\mu:V_n\cup W_n\to G$ 
by $\mu(v)=\tilde\xi(v)^{-1}$ for $v\in V_n$ and 
$\mu(w)=\tilde\eta(w)^{-1}$ for $w\in W_n$, 
we can regard ${\cal M}(V_n,W_n)/({\cal V}_n+{\cal W}_n)$ 
as the abelian group ${\cal N}(V_n,W_n)$ which was defined in Section 2. 
Let $\tilde\pi_n$ denote the connecting homomorphism 
from 
\[ {\cal N}(V_n,W_n)={\cal M}(V_n,W_n)/({\cal V}_n+{\cal W}_n) \]
to 
\[ {\cal N}(V_{n+1},W_{n+1})
={\cal M}(V_{n+1},W_{n+1})/({\cal V}_{n+1}+{\cal W}_{n+1}) \]
induced by $\pi_n:{\cal M}(V_n,W_n)\to{\cal M}(V_{n+1},W_{n+1})$. 

From the lemma above, 
the group of coinvariants $\Coinv(Z,\omega)$ is isomorphic to 
the inductive limit of $C_n/D_n$. 
Consequently, we get the following. 

\begin{lem}\label{Coinv=lim}
The group of coinvariants $\Coinv(Z,\omega)$ is isomorphic to 
the inductive limit of 
${\cal N}(V_n,W_n)={\cal M}(V_n,W_n)/({\cal V}_n+{\cal W}_n)$. 
\end{lem}

Now we are ready to prove our main theorem. 

\begin{proof}[Proof of Theorem \ref{main}]
Theorem \ref{T_G=wedge} tells us that 
the torsion subgroup $T({\cal N}(V_n,W_n))$ is 
isomorphic to $T_G\cong G\wedge G$. 
By Lemma \ref{connect}, 
the homomorphism $\pi_n$ 
from ${\cal M}(V_n,W_n)$ to ${\cal M}(V_{n+1},W_{n+1})$ is 
exactly the same as the homomorphism considered in Lemma \ref{torsioniso}. 
It follows that $\tilde\pi_n$ induces an isomorphism 
between the torsion subgroups of 
${\cal N}(V_n,W_n)$ and ${\cal N}(V_{n+1},W_{n+1})$. 
Hence, by the lemma above, 
we can conclude that 
the torsion subgroup of $\Coinv(Z,\omega)$ is isomorphic to $G\wedge G$. 
\end{proof}

\section{An example}

In \cite[Section 6.3]{GHK}, 
it was observed that 
the octagonal tiling described by the data $(V,\Gamma_8,{\cal W}^b_8)$ 
has a torsion component $\Z_2$ in the top-dimensional cohomology. 
We would like to see that 
this torsion can be explained in our context. 

Instead of $\Gamma_8$, 
it is convenient to use the lattice $\Gamma$ 
which is rotated by $\pi/8$ with respect to $\Gamma_8$, that is, 
$\Gamma$ is the $\Z$-span of the vectors 
\[ e_k=\left(\cos\frac{(2k-1)\pi}{8},\sin\frac{(2k-1)\pi}{8}\right) \]
in $\R^2$ for $k=0,1,2,3$. 
The set of representative singular lines ${\cal W}^b_8$ 
consisted of four lines. 
But, we `cut' the plane along only two directions, 
horizontal lines and vertical lines, 
because they are enough to obtain the torsion component. 
More precisely, 
letting $W_1$ (resp. $W_2$) be the horizontal (resp. vertical) line 
which passes through the origin, 
we define a totally disconnected (noncompact) space $\Pi$ 
by `cutting' the plane 
along the lines $W_i+\gamma$ for $i=1,2$ and $\gamma\in\Gamma$. 
Then $\Gamma$ acts on $\Pi$ naturally by translation. 

Let $\Gamma_0\subset\Gamma$ be the subgroup 
generated by $e_0+e_1$ and $e_2+e_3$ and 
let $Z$ be the quotient of $\Pi$ by $\Gamma_0$. 
We can identify $Z$ with the square surrounded by 
$W_i$ and $W_i+e_0+e_1+e_2+e_3$ for $i=1,2$. 
Clearly $Z$ is a Cantor set and 
the translations by $e_1$ and $e_2$ induce 
a $\Z^2$-action $\omega$ on $Z$. 
It is well-known that 
$(\Pi,\Gamma)$ has the same dynamical cohomology as $(Z,\omega)$. 

Let us look at $(Z,\omega)$. 
It is not hard to see that 
the translations by $(e_0+e_1)/2$ and $(e_2+e_3)/2$ 
are well-defined on $\Pi$. 
It follows that 
they induce a $\Z_2\times\Z_2$-action on $Z$. 
We can easily see that 
the square surrounded by $W_i$ and $W_i+(e_0+e_1+e_2+e_3)/2$ for $i=1,2$ 
is a fundamental domain for this $\Z_2\times\Z_2$-action. 
Moreover, the translations by $e_1$ and $e_2$ can be regarded 
as the translations by $e_1-(e_0+e_1)/2$ and $e_2-(e_2+e_3)/2$ 
on this small square. 
We can observe that 
$e_1-(e_0+e_1)/2$ is vertical and $e_2-(e_2+e_3)/2$ is horizontal. 
Therefore the quotient of $Z$ by the $\Z_2\times\Z_2$-action is 
the product of two Cantor minimal $\Z$-systems and 
each of them is conjugate to the same Cantor minimal $\Z$-system $\xp$. 
Here, the Cantor set $X$ is obtained 
by `cutting' the torus $\R/\Z$ at $\Z+\Z\sqrt{2}$, and 
the minimal homeomorphism $\phi$ is the translation by $\sqrt{2}$ on $X$. 
Since $(Z,\omega)$ is conjugate to 
a skew product extension of two copies of $\xp$ 
associated with $\Z_2\times\Z_2$-valued cocycles, 
by Theorem \ref{main}, 
the torsion subgroup of $\Coinv(Z,\omega)$ is isomorphic to $\Z_2$. 

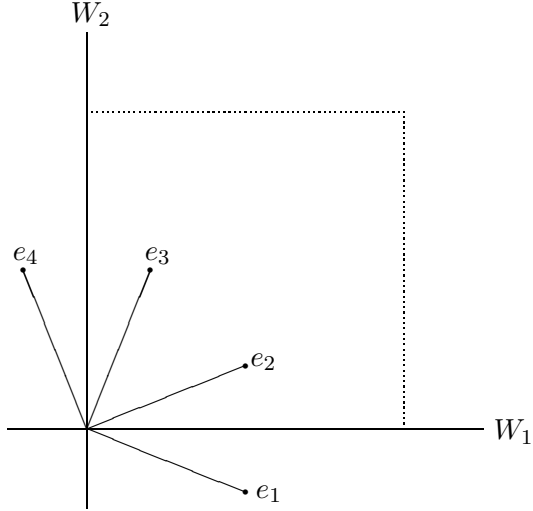
\begin{figure}
\begin{center}
\begin{picture}(200,200)

\put(30,0){\line(0,1){180}}
\put(0,30){\line(1,0){180}}
\put(24,184){$W_2$}
\put(184,26){$W_1$}

\qbezier[60](30,150)(90,150)(150,150)
\qbezier[60](150,30)(150,90)(150,150)

\put(30,30){\line(5,-2){60}}
\put(30,30){\line(5,2){60}}
\put(30,30){\line(2,5){24}}
\put(30,30){\line(-2,5){24}}

\put(90,6){\circle*{2}}
\put(90,54){\circle*{2}}
\put(54,90){\circle*{2}}
\put(6,90){\circle*{2}}

\put(94,4){$e_1$}
\put(92,54){$e_2$}
\put(52,94){$e_3$}
\put(2,94){$e_4$}

\end{picture}
\caption{The lattice $\Gamma$ and the square $Z$}
\end{center}
\end{figure}

\bigskip

\noindent
\textbf{Acknowledgment.}
I would like to thank Johannes Kellendonk and Ian Putnam 
for many helpful discussions. 
This research was carried out 
while I was visiting the University of Victoria 
under the support of the Japanese Society for the Promotion of Science.


\begin{thebibliography}{ABCD}
\bibitem[FHK]{FHK}Forrest, A.H.; Hunton, J. R.; Kellendonk, J.; 
\textit{Topological invariants for projection method patterns}, 
Mem. Amer. Math. Soc. 159 (2002), no. 758. 
\bibitem[G]{G}G\"ahler, F; 
Lecture given at PIMS Workshop on 
Aperiodic Order: Dynamical Systems, Combinatorics, and Operators, 
Banff, June 2004. 
\bibitem[GHK]{GHK}G\"ahler, F; Hunton, J. R.; Kellendonk, J.; 
\textit{Torsion in tiling homology and cohomology}, 
preprint. math-ph/0505048. 
\bibitem[HPS]{HPS}Herman, R. H.; Putnam, I. F.; Skau, C. F.; 
\textit{Ordered Bratteli diagrams, dimension groups 
and topological dynamics}, 
Internat. J. Math. 3 (1992), 827--864. 
\bibitem[M]{M}Matui, H.; 
\textit{Finite order automorphisms and dimension groups 
of Cantor minimal systems}, 
J. Math. Soc. Japan 54 (2002), 135--160. 
\end{thebibliography}
\end{document}